
\magnification=\magstep 1
\baselineskip=15pt     
\parskip=3pt plus1pt minus.5pt
\overfullrule=0pt
\font\hd=cmbx10 scaled\magstep1

\def\Pone{{\bf P}^1}
\def\P{{\bf P}}

\def\SS{{\cal S}}
\def\X{{\cal X}}
\def\O{{\cal O}}

\def\PP{{\cal P}}
\def\T{{\cal T}}

\def\G{{\cal G}}

\def\F{{\cal F}}
\def\I{{\cal I}}
\def\M{{\cal M}}

\def\E{{\cal E}}

\def\OP1{\O_{\Pone}}
\def\Pic{\mathop{\rm Pic}}

\def\im{\mathop{\rm im}}

\def\hom{\mathop{\rm Hom}}

\def\exact#1#2#3{0\rightarrow#1\rightarrow#2\rightarrow#3\rightarrow0}

\def\mapright#1{\smash{
  \mathop{\longrightarrow}\limits^{#1}}}

\def\mapdown#1{\Big\downarrow
   \rlap{$\vcenter{\hbox{$\scriptstyle#1$}}$}}
\centerline{\hd The Deformation Space of Calabi-Yau $n$-folds with}
\centerline{\hd Canonical
Singularities Can Be Obstructed.}
\medskip
\centerline{\it Mark Gross}
\medskip
\centerline{February, 1994}
\medskip
\centerline{Department of Mathematics}
\centerline{Cornell University}
\centerline{Ithaca, NY 14850}
\centerline{mgross@math.cornell.edu}
\bigskip
\bigskip
{\hd \S 0. Introduction.}

The Bogomolov-Tian-Todorov theorem ([10] and [12])
states that a non-singular $n$-fold
$X$ with $c_1(X)=0$ has unobstructed deformation theory, i.e. the moduli
space of $X$ is non-singular. This theorem was reproven
using algebraic methods by Ran in [7]. Since then, it has been proven
for Calabi-Yau $n$-folds with various mild forms of isolated singularities:
ordinary double points by Kawamata [5] and Tian [11], Kleinian singularities
by Ran [8], and finally, in the case of threefolds, arbitrary
terminal singularities by Namikawa in [6]. Now the most natural
class of singularities in the context
of Calabi-Yau
$n$-folds are canonical singularities. Indeed, if $X$ is a Calabi-Yau
$n$-fold with terminal singularities, and $f:X\rightarrow Y$ is a
birational contraction, $Y$ normal, then $Y$ has canonical singularities.
Thus the natural question to ask is: is the deformation space of Calabi-Yau
$n$-folds with canonical singularities unobstructed?

Given the history of this problem presented above, it appears
worthwhile to give a counterexample to this most general
question. We give an example of a Calabi-Yau
$n$-fold $X$ with the simplest sort of dimension 1 canonical singularities,
and show that $X$ lies in the intersection of two distinct families
of Calabi-Yau $n$-folds. One is a family of generically
non-singular Calabi-Yaus,
and the other is a family of Calabi-Yaus which generically
have terminal singularities. (In the
case $n=3$, these are also non-singular.) In particular, the point
of the moduli space corresponding to $X$ is in the intersection of two
components of moduli space, and hence has obstructed deformation theory.

We do not address the issue of isolated singularities here. That issue is
more of a local one, and the obstructedness of Calabi-Yaus with isolated
singularities is related to the obstructedness of the singularities
themselves. We will explore this in a future paper, and give applications
to smoothing Calabi-Yaus with canonical singularities.

{\it Notation:} If $\F$ is a vector bundle on a variety, we use
Grothendieck's convention for ${\bf P}(\F)$, so $\P(\F)$
denotes ${\bf Proj}(S(\F))$.
By a Calabi-Yau $n$-fold,
we mean a normal $n$-dimensional projective variety $X$ over the complex
numbers with
at worst canonical singularities, $\omega_X\cong\O_X$, and $h^1(\O_X)=0$.

{\hd \S 1. Two Families of Calabi-Yau $n$-folds.}

Let's start by defining
two distinct families of Calabi-Yau $n$-folds. Let $n$ be
an integer which is at least 3.
Let $P_1=\P(\O_{\Pone}^{\oplus (n+1)})\cong\Pone\times \P^n$, and
let
$P_2=\P(\E)$, where $\E=
\O_{\Pone}(-1)\oplus\O_{\Pone}^{\oplus (n-1)}\oplus\O_{\Pone}(1)$.
For each $P_i$, the Picard group is generated by $t$, the class of
$\O_{P_i}(1)$, and $f$, the class of a fibre of the projection $\pi:P_i
\rightarrow\Pone$. The canonical class of $P_i$ is then $K_{P_i}=-(n+1)t-2f$.
For an element $s\in H^0(\omega_{P_i}^{-1})$, we denote the
zero locus of $s$ by $X_i(s)\subseteq
P_i$.
For a general $s$ we have $X_1(s)$ non-singular. However, for
$n>3$, this is
not true for $X_2(s)$.

We need to examine the structure of the singularities of $X_2(s)$. First,
let's look at $P_2$ in more detail.
There is a section $C\subseteq P_2$
of the bundle $\pi:P_2\rightarrow \Pone$ given by the inclusion
$\P(\O_{\Pone}(-1))\subseteq \P(\E)$ induced by the quotient
map $\E\rightarrow \O_{\Pone}(-1)\rightarrow 0$. We denote by $\I_C^p$
the $p$th power of the ideal sheaf of $C$ in $P_2$.

Let $F^p\subseteq S^{n+1}\E$ be the subbundle given by
$$F^p=\bigoplus_{i=p}^{n+1} S^i(\OP1(-1))\otimes S^{n+1-i}(\OP1^{n-1}
\oplus\OP1(1)).$$
This yields a filtration of $S^{n+1}\E$:
$$S^{n+1}\E=F^0\supseteq F^1\supseteq\cdots\supseteq F^{n+1}
\supseteq F^{n+2}=0$$
such that $F^p/F^{p+1}\cong S^p(\OP1(-1))\otimes S^{n+1-p}(\OP1^{n-1}
\oplus \OP1(1))$. This is the natural filtration on
$S^{n+1}\E$ induced by the exact sequence
$$\exact{\OP1(-1)}{\E}{\OP1^{n-1}\oplus\OP1(1)}$$
via [3], II, Ex. 5.16 c).

Let
$$G^p=\bigoplus_{i=0}^{n+1-p} S^i(\OP1(-1))\otimes S^{n+1-i}(\OP1^{n-1}
\oplus\OP1(1)).$$
We have
$S^{n+1}\E\cong G^{n+2-p}\oplus F^p$.

\proclaim Lemma 1.1.
\item{(i)} There is a natural isomorphism
$\pi_*(\I_C^p\otimes\omega_{P_2}^{-1})
\cong G^p(2)$ inducing a commuting diagram
$$\matrix{H^0(\I_C^p\otimes\omega_{P_2}^{-1})&\subseteq&H^0(\omega_{P_2}^{-1})
\cr
\mapdown{\cong}&&
\mapdown{\cong}\cr
H^0(G^p(2))&\subseteq&H^0((S^{n+1}\E)(2))\cr}$$
\item{(ii)}
Let $V=H^0(\I_C^{n-1}\otimes \omega_{P_2}^{-1})
\subseteq H^0(\omega_{P_2}^{-1})$.
If $s\in V$ is a general element, then $X_2(s)$ is a Calabi-Yau
$n$-fold with canonical
singularities along $C$ and is non-singular elsewhere. In addition,
the natural map
$$\psi:S^{n-1}H^0(\O_{P_2}(t))\otimes S^2H^0(\O_{P_2}(t+f))
\rightarrow H^0(\omega_{P_2}^{-1})$$
has image $\im \psi=V\subseteq H^0(\omega_{P_2}^{-1})$.
\item{(iii)} If
$s\in H^0(\omega^{-1}_{P_2})$ is a general element, then $X_2(s)$
is non-singular outside of $C$ and has singularities generically of
multiplicity $\lfloor {n\over 2} \rfloor$ along $C$.
If $n=3$, $X_2(s)$ will be non-singular.

Proof: (i) Let $V_m$ be an $m$-dimensional vector space, and let
$P\in \P(V)$ be a point. Giving $P$ is the same thing as giving a
one-dimensional quotient space $V_1$ of $V_m$, or an exact sequence
$$\exact{V_{m-1}}{V_{m}}{V_1}.\leqno{(*)}$$
Now $V_m=H^0(\O_{\P(V_m)}(1))$, and $V_{m-1}\subseteq V_m$ is the
subset of linear forms which vanish at the point $P$.
We have  a filtration of $S^dV_m=H^0(\O_{\P(V_m)}(d))$,
$$S^dV_m=W^0\supseteq\cdots\supseteq W^{d+1}=0,$$
with $W^p/W^{p+1}\cong S^pV_{m-1}\otimes S^{d-p}V_1$.
$W^p$ in fact then consists of $d$-forms vanishing to order
at least $p$ at $P$, so we see that $W^p=H^0(\I_P^p(d))$.

This can all be relativized. In the situation of part (i), the
curve $C$ comes from the split exact sequence
$$\exact{\OP1(1)\oplus\OP1^{n-1}}{\E}{\OP1(-1)},\leqno{(**)}$$
which corresponds to the sequence (*).
The corresponding filtration of $S^{n+1}\E$ is then
$$S^{n+1}\E=G^0\supseteq\cdots\supseteq G^{n+2}=0.$$
The relativized statement corresponding to $H^0(\I_P^p(d))=W^p$
is $\pi_*\I_C^p((n+1)t)\cong G^p$, and so $\pi_*\I_C^p\otimes\omega_{P_2}^{-1}
\cong G^p(2)$.
Now $H^0(\omega_{P_2}^{-1})\cong H^0(\pi_*\omega_{P_2}^{-1})
=H^0((S^{n+1}\E)(2))$, and $H^0(\I_C^p\otimes\omega_{P_2}^{-1})
\cong H^0(\pi_*(\I_C^p\otimes\omega_{P_2}^{-1}))
\cong H^0(G_p(2))\subseteq H^0((S^{n+1}\E)(2))$. This gives the desired
diagram.

(ii) Let $s$ be a general element of
$V$. For any $s$, we have $H^1(\O_{X_2(s)})\cong H^2(\omega_{P_2})=0$.
Thus, to show that $X_2(s)$ is a Calabi-Yau $n$-fold
with canonical singularities, since $K_{X_2(s)}=0$ by adjunction,
it is enough to show that there is
a resolution of singularities $\tilde X_2(s)\rightarrow X_2(s)$
such that $K_{\tilde X_2(s)}=0$.

Let $b:\tilde P_2\rightarrow P_2$ be the blow-up of $P_2$ along $C$
with exceptional divisor $E$. The proper transform $\tilde X_2(s)$ of
$X_2(s)$ for $s\in V$ will be an element of the linear system
$|(n+1)b^*t+2b^*f-(n-1)E|$. If this linear system is base-point-free,
then for general $s\in V$, $\tilde X_2(s)\rightarrow X_2(s)$ will
be a resolution of singularities. Furthermore, $K_{\tilde P_2}
=-(n+1)b^*t-2b^*f+(n-1)E$, so $K_{\tilde X_2(s)}=0$.

To show that $|(n+1)b^*t+2b^*f-(n-1)E|$ is
base-point-free, it is enough to show that
$|b^*t-E|$ is base-point-free and that
$|t+f|$ is base-point-free, for then so is
$|(n-1)(b^*t-E)+2b^*(t+f)|$.

It is easy to see, in general, that if
$$\E=\bigoplus_{i=1}^n\O_{\Pone}(a_i)$$ is a vector bundle
over $\Pone$, and $t=c_1(\O_{\P(\E)}(1))$, then the base locus
of $|t|$ is $\P(\F)\subseteq \P(\E)$, where
$$\F=\bigoplus_{\hbox{$i$ with $a_i<0$}} \O_{\Pone}(a_i),$$
with the inclusion induced
by the natural surjection $\E\rightarrow \F$. Thus, in particular,
for $P_2$,
we see that the base locus of $|t+f|$ is empty, and the base locus
of $|t|$ is the curve $C$.

Now to see that $|b^*t-E|$ is base-point-free, note that
we have an exact sequence
$$\exact{\O_{P_2}(t-f)}{\O_{P_2}(t)}{\O_f(t)}$$
with $f\cong \P^n$.  Now $h^1(\O_{P_2}(t-f))=h^1(\pi_*\O_{P_2}(t-f))
=h^1(\E(-1))=1$, and similarly $h^1(\O_{P_2}(t))=0$, so the image of
$H^0(\O_{P_2}(t))$ in $H^0(\O_{f}(t))=H^0(\O_{\P^n}(1))$ is codimension
one. This image must yield the linear system of hyperplanes in $f$ which
contain the point $C\cap f$. Now after blowing up $C$,
the linear system $|b^*t-E|$ is isomorphic, via proper transform,
to $|t|$, and thus its restriction to the proper transform of $f$
is now base-point-free. Since there is a divisor of type $f$ through
any point in $P_2$, this shows that $|b^*t-E|$ is base-point-free.
This proves the first statement.

For the second, first note that the linear system $|\im\psi|\subseteq
|-K_{P_2}|$ is spanned by reducible divisors consisting of a union of
$n-1$ divisors in $|t|$ and 2 divisors in $|t+f|$. We have seen that
$C$ is contained in any divisor in $|t|$, so that this reducible
divisor in $|-K_{P_2}|$ contains $C$ to order at least $n-1$.
Thus $\im\psi\subseteq V$.

Let $\E'=\OP1^{n-1}\oplus\OP1(1)$. Then
$$G^{n-1}(2)\cong (S^{n+1}\E')(2)\oplus (S^n\E')(1)\oplus S^{n-1}(\E'),$$
and since $\E'$ is generated by global sections over $\Pone$, it is easy to see
that the maps
$$\eqalign{V_1=S^{n+1}(H^0(\E'))\otimes H^0(\OP1(2))&\rightarrow
H^0((S^{n+1}\E')(2))\subseteq V,\cr
V_2=S^n(H^0(\E'))\otimes H^0(\OP1(1))&\rightarrow
H^0((S^n\E')(1))\subseteq V,\cr}$$
and
$$V_3=S^{n-1}(H^0(\E'))\otimes H^0(\OP1)\rightarrow
H^0(S^{n-1}\E')\subseteq V$$
are all surjective, and so the map $\psi':V_1\oplus V_2\oplus V_3\rightarrow
V$ is surjective.

Now the linear system $|\psi'(V_1)|\subseteq |V|$ is spanned by
divisors consisting of a sum of $n+1$ divisors in $|t|$ and two in
$|f|$; the linear system $|\psi'(V_2)|\subseteq |V|$ is spanned by divisors
consisting of a sum of $n$ divisors in $|t|$, one in $|f|$, and the
divisor $D\in |t+f|$ given by the inclusion of $\OP1(-1)$ in $\E$ coming from
the splitting of the sequence (**). The linear system $|\psi'(V_3)|\subseteq
|V|$ is spanned by divisors consisting of a sum of $n-1$ divisors in $|t|$
and $2D$. All these divisors are contained in $|\im\psi|$, so $V\subseteq
\im \psi$. Thus $V=\im \psi$.

(iii) By (i), we  have a filtration
$$H^0(\omega_{P_2}^{-1})=H^0(G^0(2))\supseteq H^0(G^1(2))\supseteq
\cdots\supseteq H^0(G^{n+2}(2))=0,$$
with $H^0(G^p(2))$ consisting of those sections of $H^0(\omega_{P_2}^{-1})$
which vanish to order at least $p$ along $C$.
A simple calculation shows that
$$H^0(G^{\lfloor n/2\rfloor}(2))=
H^0(G^{\lfloor n/2\rfloor-1}(2))=\cdots=H^0(G^0(2)),$$
but that $$H^0(G^{\lfloor n/2\rfloor+1}(2))\neq
H^0(G^{\lfloor n/2 \rfloor}(2)).$$
Thus the general element of $|-K_{P_2}|$ has singularities
generically of multiplicity $\lfloor n/2\rfloor$ along $C$. Part (ii)
shows that the general element of $|-K_{P_2}|$ has no
singularities outside of $C$.

We now restrict to the case $n=3$. Let $S\subseteq P_2$, $S\cong
\P(\OP1(-1)\oplus\OP1(1)),$ be the surface determined by a surjection
$\E\rightarrow \OP1(-1)\oplus\OP1(1)=\E''$. As the exact sequence
$$\exact{\OP1^2}{\E}{\E''}$$
is split, $S^4\E''$ is a direct summand of $S^4\E$. The map
$H^0(S^4\E\otimes\OP1(2))\rightarrow H^0(S^4\E''\otimes\OP1(2))$
is then surjective, and this coincides with the restriction map
$H^0(\omega_{P_2}^{-1})\rightarrow H^0(\omega_{P_2}^{-1}|_S)$
via the diagram
$$\matrix{H^0(\omega_{P_2}^{-1})&\cong&H^0(\pi_*\omega_{P_2}^{-1})&\cong&
H^0(S^4\E\otimes\O_{\Pone}(2))\cr
\mapdown{}&&\mapdown{}&&\mapdown{}\cr
H^0(\omega_{P_2}^{-1}|_{S})&\cong&H^0(\pi_*\omega_{P_2}^{-1}|_{S})&\cong&
H^0(S^4\E''\otimes\O_{\Pone}(2))\cr}$$
Now $-K_{P_2}|_S\sim 4C+6f$ where $\Pic S$ is spanned by
$C$, which has self-intersection $-2$ on $S$, and $f$.
By [3],
V, 2.18, a general member of $|4C+6f|$ consists of a sum of $C$
and an irreducible non-singular curve of type $3C+6f$, disjoint from
$C$. Thus the general element of $|4C+6f|$ is non-singular,
and $X_2(s)\cap S$ is non-singular for general $s\in H^0(\omega_{P_2}^{-1})$.
Thus $X_2(s)$ is non-singular along $C$.
$\bullet$

{\hd \S 2. The Example.}

We are now ready to give our example of a Calabi-Yau $n$-fold with
canonical singularities with singular Kuranishi space. A versal
Kuranishi space exists for any compact complex space
by [1] or [2]. By our definition of a Calabi-Yau $n$-fold and [4], (8.6),
$\hom(\Omega_X^1,\O_X)\cong H^1(\O_X)=0$, and so by [9] this versal
Kuranishi space is universal for a Calabi-Yau $n$-fold.

Returning to the setup of \S 1,
as $Ext^1(\O_{\Pone}(1)\oplus\O_{\Pone}^{n-1},
\O_{\Pone}(-1))$ is one dimensional,
there is a universal extension bundle $\F$ on ${\bf A}^1\times\Pone$,
for which $\F|_{0\times\Pone}\cong \O_{\Pone}(-1)\oplus
\O_{\Pone}^{n-1}\oplus\O_{\Pone}(1)$, and $\F|_{t\times \Pone}=
\O_{\Pone}^{n+1}$ for $t\in {\bf A}^1-\{0\}$.
This yields a family of $\P^n$-bundles over $\Pone$ via $\P(\F)
\rightarrow {\bf A}^1\times\Pone$. If we take
$\X\subseteq\P(\F)$ to be the zero locus of a section of
$\omega^{-1}_{\P(\F)/{\bf A}^1}$, we would presumably obtain a family
of Calabi-Yau $n$-folds $\X\rightarrow {\bf A}^1$, the general fibre
being contained in $P_1$, but the fibre over $0\in
{\bf A}^1$ being contained in $P_2$. Now we can apply the following
Lemma.

\proclaim Lemma 2.1. If there exists a flat family $\X\rightarrow \SS$ with
$0\in \SS$
a point with $\X_0\cong X_2(s)$ for some $s\in H^0(\omega_{P_2}^{-1})$ with
$\X_0$ a Calabi-Yau $n$-fold,
and $\X_t$ isomorphic to a non-singular member of $|-K_{P_1}|$
for $t\in \SS-\{0\}$, then the Kuranishi space at $\X_0$ is
singular.

Proof: The Kuranishi space of $\X_0$ must contain a subspace $\M_1$
corresponding to deformations of $\X_0$ to non-singular elements
of $|-K_{P_1}|$, and a subspace $\M_2$ corresponding to deformations of
$\X_0$ to elements of $|-K_{P_2}|$.

$X_1(s')$ is non-singular for general $s'\in H^0(\omega_{P_1}^{-1})$,
so the dimension
of the Kuranishi space  for $X_1(s')$ can be calculated by calculating
the dimension of its tangent space, which is $H^1(\T_{X_1(s')})$. A simple
calculation shows that this coincides with the value for the dimension
of the Kuranishi space one would expect via a naive dimension counting
of the number of moduli in $P_1$:
$$\dim\M_1=h^1(\T_{X_1(s)})=h^0(\omega_{P_1}^{-1})-1-\dim Aut(P_1).$$
Furthermore $\M_1$ must be an irreducible component of the Kuranishi
space of $\X_0$.

Now the dimension of $\M_2$ is at least
$$h^0(\omega_{P_2}^{-1})-1-\dim Aut(P_2).$$
An automorphism of $P_2$ is induced by an automorphism of $\Pone$ and
an automorphism of the bundle $\E$. An automorphism of $\E$ is induced
by an $(n+1)\times (n+1)$ matrix of forms over $\Pone$, of which
$(n-1)^2+2$ of these entries are constant forms, $2(n-1)$ are linear
forms, and one is a quadratic form. The other entries must be zero.
This gives a dimension $(n+1)^2+1$ set of matrices, or $(n+1)^2$
dimensional modulo scalars.
Thus $\dim Aut(P_2)=(n+1)^2+3=\dim Aut(P_1)+1$.
Meanwhile $h^0(\omega_{P_2}^{-1})\ge h^0(\omega_{P_1}^{-1})+1,$
with equality holding if and only if $n=3$.
So, for $n>3$, $\dim\M_2>\dim\M_1$, and the Kuranishi space must
have at least two irreducible components meeting at $\X_0$, $\M_1$
being one of them, and the other containing $\M_2$.

If $n=3$, $\M_1$ and $\M_2$ are the same dimension, so this
argument does not suffice. However, if $\M_1$ and $\M_2$ coincide
then there would be a non-singular Calabi-Yau $X_2(s_2)$ isomorphic
to $X_1(s_1)$ for some $s_1$, $s_2$ via an isomorphism $\alpha:
X_1(s_1)\rightarrow X_2(s_2)$. But such an isomorphism would have to
preserve the cubic intersection form, with $(\alpha^* D)^3=D^3$
for $D\in\Pic X_2(s_2)$. It is then easy to see that $\alpha^*t=t$
and $\alpha^*f=f$ is the only possibility for $\alpha^*:\Pic X_2(s_2)
\rightarrow \Pic X_1(s_1)$. But $t$ is a nef divisor on $X_1(s_1)$
but not on $X_2(s_2)$, so there is no such isomorphism. Thus
$\M_1$ and $\M_2$ are two distinct components of the Kuranishi
space at $\X_0$.

Thus, in any event, the Kuranishi space is reducible,
hence singular, at $\X_0$.
$\bullet$

Thus, to construct our desired counterexample, we just need to get
control of the singularities of $\X_0$ to ensure that they
are no worse than canonical singularities. We do this by showing that we
can construct a family $\X\rightarrow{\bf A}^1$ as above with $\X_0
\cong X_2(s)$ for any $s\in V$.

\proclaim Theorem 2.2. Let $s\in V\subseteq H^0(\omega_{P_2}^{-1})$.
Then the Kuranishi space of $X_2(s)$ is singular at $X_2(s)$.

Proof. By Lemma 2.1, it is enough to show that $X_2(s)$ is deformation
equivalent to $X_1(s')$ for general $s'$. We use the following
construction:

Let $\F$ be the universal extension bundle over $T={\bf A}^1\times\Pone$,
as at the beginning of this section.
Let $p_1$ and $p_2$ be the projections of $T$ onto the first and
second factors respectively.  We set $\O_T(1)=p_2^*\OP1(1)$.
Let $\PP=\P(\F)$, $\pi:\PP\rightarrow T$ the projection, and
denote by $t$ the class of $\O_{\PP}(1)$, and $f$ the class of $\pi^*\O_T(1)$.

\proclaim Lemma 2.3. Let $x\in{\bf A}^1$, and let
$$\phi_x:p_{1*}\pi_*\O_{\PP}((n+1)t+2f)\otimes k(x)
\rightarrow H^0(\pi_*\O_{\PP}((n+1)t+2f)|_{x\times\Pone})$$
be the restriction map, where $k(x)$ is the
residue field of ${\bf A}^1$ at $x$.
Then $\phi_x$ is surjective if $x\not=0$,
and $\im \phi_0=V\subseteq H^0(\omega_{P_2}^{-1})
=H^0(\pi_*\O_{\PP}((n+1)t+2f)|_{0\times\Pone}),$
where $V$ is the subspace of Lemma 1.1, (ii).

Proof: First we have
$\pi_*\O_{\PP}((n+1)t+2f)=(S^{n+1}\F)(2)$.
For $x\not=0$, $(S^{n+1}\F)(2)|_{x\times\Pone}
=\O_{\Pone}(2)^N$ for suitable $N$, so by Grauert's Theorem,
[3], III, 12.9,
$\phi_x$ is surjective.

For $x=0$, first note that the maps
$$p_{1*}\pi_*\O_{\PP}(t)\otimes k(x)
\rightarrow H^0(\pi_*\O_{\PP}(t)|_{0\times\Pone})$$
and
$$p_{1*}\pi_*\O_{\PP}(t+f)\otimes k(x)
\rightarrow H^0(\pi_*\O_{\PP}(t+f)|_{0\times\Pone})$$
are surjective, again by Grauert's Theorem. By Lemma 1.1, (ii),
this shows that $\im \phi_0$ contains $V$.

By construction, we have
$$\exact{\O_T(-1)}{\F}{\O_T^{n-1}\oplus\O_T(1)}.$$
Let
$$S^{n+1}\F\supseteq \F^0\supseteq\cdots\supseteq\F^{n+2}=0$$
be the filtration of $S^{n+1}\F$ induced by the above extension,
and let
$$\G^p=S^{n+1}\F/\F^{n+2-p}.$$

For $p=n-1$, this yields the sequence
$$\exact{\F^3(2)}{(S^{n+1}\F)(2)}{\G^{n-1}(2)}.$$
Restricting this to $0\times\Pone$, this sequence splits to
obtain
$$S^{n+1}\E\cong F^3(2)\oplus G^{n-1}(2).$$

Now $p_{1*}\F^3(2)=0$ since $\F^3(2)|_{x\times\Pone}\cong \OP1(-1)^N$ for
some suitable $N$,
$x\not=0$, and the map
$H^0(\G^{n-1}(2))\otimes k(0)\rightarrow H^0(G^{n-1}(2))$ is an
isomorphism by Grauert's Theorem. This yields the following diagram:
$$\matrix{&&0&\mapright{}
&p_{1*}(S^{n+1}\F)(2)\otimes k(0)&\mapright{\alpha}
&p_{1*}\G^{n-2}(2)\otimes k(0)\cr
&&&&\mapdown{\phi_0}&&\mapdown{\beta\cong}\cr
0&\mapright{}&H^0(F^3(2))&\mapright{}&H^0((S^{n+1}\E)(2))&
\mapright{}&H^0(G^{n-1}(2))\cr}$$
Since $\im\phi_0\supseteq V=H^0(G^{n-1}(2))$, yet at the same
time this diagram shows $\phi_0$ injects into $H^0(G^{n-1}(2))$,
we see $\phi_0$ has as its image exactly $V=H^0(G^{n-1}(2))$.
$\bullet$

To prove Theorem 2.2,
let $s_0\in V\subseteq H^0((S^{n+1}\F)(2)|_{0\times
\Pone})$ be any element. We can then lift $\phi_0^{-1}(s_0)$ to
a section $s$ of $S^{n+1}(\F)(2)$
such that $s|_{1\times \Pone}=s_1\in H^0((S^{n+1}\F)(2)|_{1\times\Pone})$
is a general section, and $s|_{0\times\Pone}=s_0$, by Lemma
2.3. Let $\X\subseteq\PP$ be the zero-locus of the corresponding
section of $\O_{\PP}((n+1)t+2f)$. The projection $\X\rightarrow {\bf A}^1$
gives a family of Calabi-Yaus over ${\bf A}^1$, with $\X_1=X_1(s_1)$,
and $\X_0=X_2(s_0)$. $\bullet$

{\hd Bibliography}

\item{[1]} Douady, A., ``Le Probl\`eme des Modules Locaux pour les
Espaces ${\bf C}$-Analytiques Compacts,'' {\it Ann. scient. \'Ec. Norm.
Sup.,} 4e s\'erie, {\bf 7}, (1974) 569--602.
\item{[2]} Grauert, H., ``Der Satz von Kuranishi f\"ur Kompakte Komplexe
R\"aume,'' {\it Inv. Math.,} {\bf 25}, (1974) 107--142.
\item{[3]} Hartshorne, R., {\it Algebraic Geometry,} Springer-Verlag 1977.
\item{[4]} Kawamata, Y., ``Minimal Models and the Kodaira Dimension
of Algebraic Fiber Spaces,'' {\it J. Reine Angew. Math,} {\bf 363}, (1985)
1--46.
\item{[5]} Kawamata, Y., ``Unobstructed Deformations--- A Remark on a Paper
of Z. Ran,''
{\it J. Algebraic Geometry,} {\bf 1}, (1992) 183--190.
\item{[6]} Namikawa, Y., ``On Deformations of Calabi-Yau Threefolds with
Terminal
Singularities,'' preprint, 1993.
\item{[7]} Ran, Z., ``Deformations of Manifolds with Torsion or Negative
Canonical
Bundle,'' {\it J. of Algebraic Geometry,} {\bf 1}, (1992) 279--291.
\item{[8]} Ran, Z., ``Deformations of Calabi-Yau Kleinfolds,''
in {\it Essays
in Mirror Symmetry,} (ed. S.-T. Yau) Int. Press, Hong Kong, (1992) 451--457.
\item{[9]} Schlessinger, M., ``Functors on Artin Rings,'' {\it Trans.
Amer. Math. Soc.,} {\bf 130}, (1968), 208-222.
\item{[10]} Tian, G., ``Smoothness of the Universal Deformation Space of
Compact
Calabi-Yau Manifolds and its Petersson-Weil Metric,'' in {\it Mathematical
Aspects of String Theory,} 629-646, ed. S.-T. Yau, World Scientific, Singapore,
1987.
\item{[11]} Tian, G., ``Smoothing 3-folds with trivial canonical bundle
and ordinary double points,''
in {\it Essays
in Mirror Symmetry,} (ed. S.-T. Yau) Int. Press, Hong Kong, (1992) 458--479.
\item{[12]} Todorov, A., ``The Weil-Petersson Geometry of the Moduli Space of
$SU(n\ge 3)$
(Calabi-Yau) Manifolds I,'' {\it Commun. Math. Phys} {\bf 126}, (1989) 325-346.
\end